\documentclass[11pt, reqno]{amsart}

\usepackage{amsmath,amsfonts,amssymb,amsthm,enumerate,hyperref,multicol}
\begin{document}
\leftline{ \scriptsize \it }

\title[Moment Estimates]{A Note on Jain basis functions}
\maketitle

\begin{center}
{\bf G. C. Greubel} \\
Newport News, VA, United States \\
jthomae@gmail.com
\end{center}

\vspace{3mm}

\noindent \textbf{Abstract.} In the present article the moments associated with the Jain basis functions 
are developed to order ten. The moments are shown to be polynomials in one variable with polynomial 
coefficients. The polynomial coefficients are investigated and shown to be linked to sequences found 
in the On-line Encyclopedia of Integer Sequences.


\smallskip
\noindent \textbf{Oeis Sequences.} A000217, A000292, A000332, A000389, A000579, A000580, A000581, \\
\indent A000582, A000914, A000915, A001287, A001303, A008292, A008517, A053567.

\section{Introduction}
G. C. Jain, \cite{Jain}, introduced the operators
\begin{align}\label{e1} 
P_{n}^{\beta}(f,x) = \sum_{k=0}^{\infty} f\left(\frac{k}{n}\right) \, L_{n,k}^{(\beta)}(x), \hspace{5mm}
x \in [0, \infty),
\end{align}
where $0 \leq \beta < 1$ and the basis functions are defined as
\begin{align}\label{e2} 
L_{n,k}^{(\beta)}(x) = \frac{n x \, (n x + k \beta)^{k-1}}{k!} \, e^{-(n x + k \beta)}, 
\end{align}
and has the property $\sum_{k=0}^{\infty} L_{n,k}^{(\beta)}(x) = 1$ are defined for $0 \leq \beta < 1$ and reduce
to the Szasz-Mirakyan operators when $\beta = 0$. Jain introduced the method to obtain the moments to be used
with these operators. In order to define the moments of the generalized Szasz-Mirakyan operators the series
$S(r,\alpha, \beta)$ is to be defined as
\begin{align}\label{e3} 
S(r,\alpha, \beta) = \sum_{k=0}^{\infty} \frac{1}{k!} \, (\alpha + \beta k)^{k + r -1} \, e^{-(\alpha + \beta k)}, 
\hspace{5mm} r = 0, 1, 2, 3, \cdots
\end{align}
where $\alpha \, S(0, \alpha, \beta) = 1$ and
\begin{align}\label{e4} 
S(r, \alpha, \beta) &= \sum_{k=0}^{\infty} \beta^{k} \, (\alpha + \beta k) \, S(r-1, \alpha + \beta k, \beta). 
\end{align}
When $f(t) = t^{m}$ in (\ref{e1}), where $P_{n}^{\beta}(f=t^{m}; x) = B_{n}^{\beta}(t^{m}, x)$, the operators become
the desired moments used in the referenced papers. The first couple of moments are obtained by following: 
$B_{n}^{\beta}(1,x) = 1$, 
\begin{align*}
B_{n}^{\beta}(t; x) = n x \, \sum_{k=1}^{\infty} \frac{(n x + k \beta)^{k-1}}{k!} \, e^{-(n x + k \beta)} \, \left(
\frac{k}{n} \right) = x \, S(1, n x + \beta, \beta)
\end{align*}
and
\begin{align*}
B_{n}^{\beta}(t^{2}; x) = n x \, \sum_{k=0}^{\infty} \frac{(n x + k \beta)^{k-1}}{k!} \, e^{-(n x + k \beta)} \,
\left( \frac{k}{n}\right)^{2} = \frac{x}{n} \, \left[ S(2, n x + 2 \beta, \beta) + S(1, n x + \beta, \beta) \right].
\end{align*}

Jain presented the first three moments which, for most cases, is enough to establish convergence of the operators 
and be used to identify relations dependent upon these moments. In recent works, Gupta \& Greubel \cite{Gupta-Greubel-1,
Gupta-Greubel-2}, Gupta \& Malik \cite{Gupta-Malik}, and others have extended the list of moments generaly used to
fifth order. Extending the list of moments helps establish higher orders of convergence and can lead to new identities
dependent upon these moments. In this article the list of moments is extended to order ten. The coefficients, of the
variable $\beta$ will be given and examined. The coefficients of these polynomials will be linked to sequences given
in the On-line Encyclopedia of Integer Sequences (Oeis) \cite{Oeis}. 

\vspace{5mm}

\section{Developement of \texorpdfstring{S(r, $\alpha$, $\beta$)}{Lg}}
First define the Eulerian polynomials, $A_{n}(x)$, by
\begin{align}\label{e5} 
\sum_{k=0}^{\infty} k^{n} \, x^{k} = \frac{x \, A_{n}(x)}{(1-x)^{n+1}}, 
\end{align}
where $A_{n}(x) = \sum_{j=0}^{n} A(n,m) \, x^{m}$ and $A(n,m)$ are the Eulerian numbers. Eulerian numbers
are listed as sequence A008292 in Oeis and the first few polynomials are given by $A_{0}(x) = 1$ and
\begin{align}\label{e6} 
A_{1}(x) &= 1 \nonumber\\
A_{2}(x) &= 1 + x \nonumber\\
A_{3}(x) &= 1 + 4 \, x + x^2 \\
A_{4}(x) &= 1 + 11 \, x + 11 \, x^2 + x^3 \nonumber\\
A_{5}(x) &= 1 + 26 \, x + 66 \, x^2 + 26 \, x^3 + x^4. \nonumber
\end{align}
Beginning with (\ref{e3}, \ref{e4}) it can be seen that 
\begin{align}\label{e7} 
S(1, \alpha, \beta) &= \sum_{k=0}^{\infty} \beta^{k} \, (\alpha + \beta k) \, S(0, \alpha + \beta k, \beta)
= \sum_{k=0}^{\infty} \beta^{k} = \frac{1}{1-\beta} = p.
\end{align}
The next couple are given by the following:
\begin{align}\label{e8} 
S(2, \alpha, \beta) &= \sum_{k=0}^{\infty} \beta^{k} \, (\alpha + \beta k) \, S(1, \alpha + \beta k, \beta) 
= p \, \sum_{k=0}^{\infty} \beta^{k} \, (\alpha + \beta k) \nonumber\\
&= p \, (\alpha \, p \, A_{0}(\beta) + \beta^{2} \, p^{2} \, A_{1}(\beta)) 
= p \, (\alpha \, p + \beta^{2} \, p^{2}),
\end{align}
\begin{align}\label{e9} 
S(3, \alpha, \beta) &= \sum_{k=0}^{\infty} \beta^{k} \, (\alpha + \beta k) \, S(2, \alpha + \beta k, \beta)
= p^{2} \, \sum_{k=0}^{\infty} \beta^{k} \, [ (\alpha + \beta k)^{2} + p^{3} \, (\alpha + \beta k) \, \beta^{2} ]
\nonumber \\
&= p^{2} \, (\alpha^{2} \, p \, A_{0}(\beta) + 2 \alpha \, \beta^{2} \, A_{1}(\beta) + p^{3} \, \beta^{3}
\, A_{2}(\beta)) + p^{3} \, \beta^{2} \, (\alpha \, p \, \beta \, A_{0}(\beta) + \beta^{2} \, p^{2} \, A_{1}(\beta))
\nonumber\\
&= \alpha^{2} \, p^{3} + 3 \, \alpha \, \beta^{2} \, p^{4} + \beta^{3} \, (1 + 2 \beta) \, p^{5}.
\end{align}

Repeating the process a list can be generated for the polynomials $S(r, \alpha, \beta)$. The polynomials are given 
by:
\begin{align}\label{e10} 
S(0, \alpha, \beta) &= \frac{1}{\alpha} \nonumber \\
S(1, \alpha, \beta) &= p = \frac{1}{1-\beta} \nonumber\\
S(r, \alpha, \beta) &= p^{r} \left[ \alpha^{r-1} + \sum_{k=1}^{r-1} \theta_{k}^{r}(\beta) \, \alpha^{r-k-1}
\, \beta^{k+1} \, p^{k} \right], \hspace{5mm} r \geq 2.
\end{align}
It is quickly determined that $\theta_{1}^{r}(\beta) \in \{ 0, 1, 3, 6, 10, 15, 21, 28, 36, 45, \cdots \}_{r \geq 1}$
which is sequence A000217 and is the binomial coefficients $\binom{r}{2}$, ie. $\theta_{1}^{r}(\beta) = \binom{r}{2}$.

\begin{align*}
\begin{tabular}{c|ccc}
\multicolumn{1}{c}{$r$} & \multicolumn{1}{c}{$\theta_{2}^{r}(\beta)$} &
\multicolumn{1}{c}{$\theta_{3}^{r}(\beta)$} & \multicolumn{1}{c}{$\theta_{4}^{r}(\beta)$} \\
\hline
2 & 0 & 0 & 0 \\
3 & 1 + 2$\beta$ & 0 & 0 \\
4 & 4 + 11 $\beta$ & 1 + 8 $\beta$ + 6 $\beta^{2}$ & 0 \\
5 & 10 + 35 $\beta$ & 5 + 50 $\beta$ + 50 $\beta^{2}$ & 1 + 22 $\beta$ + 58 $\beta^{2}$ + 24 $\beta^{3}$ \\
6 & 20 + 85 $\beta$ & 15 + 180 $\beta$ + 225 $\beta^{2}$ & 6 + 157 $\beta$ + 508 $\beta^{2}$ + 274 $\beta^{3}$ \\
7 & 35 + 175 $\beta$ & 35 + 490 $\beta$ + 735 $\beta^{2}$ & 21 + 637 $\beta$ + 2443 $\beta^{2}$ + 1624 $\beta^{3}$ \\
8 & 56 + 322 $\beta$ & 70 + 1120 $\beta$ + 1960 $\beta^{2}$ & 56 + 1932 $\beta$ + 8568 $\beta^{2}$ + 6769 $\beta^{3}$ \\
9 & 84 + 546 $\beta$ & 126 + 2268 $\beta$ + 4536 $\beta^{2}$ & 126 + 4872 $\beta$ + 24528 $\beta^{2}$ + 22449 $\beta^{3}$ \\
10 & 120 + 870 $\beta$ & 210 + 4200 $\beta$ + 9450 $\beta^{2}$ & 252 + 10794 $\beta$ + 60816 $\beta^{2}$ + 63273 $\beta^{3}$ \\
\end{tabular}
\end{align*}

\begin{align*}
\theta_{5}^{6}(\beta) &= 1 + 52 \beta + 328 \beta^{2} + 444 \beta^{3} + 120 \beta^{4} \\
\theta_{5}^{7}(\beta) &= 7 + 420 \beta + 3108 \beta^{2} + 5096 \beta^{3} + 1764 \beta^{4}  \\
\theta_{5}^{8}(\beta) &= 28 + 1904 \beta + 16170 \beta^{2} + 31136 \beta^{3} + 13132 \beta^{4}  \\
\theta_{5}^{9}(\beta) &= 84 + 6384 \beta + 61194 \beta^{2} + 135324 \beta^{3} + 67284 \beta^{4}  \\
\theta_{5}^{10}(\beta) &= 210 + 17640 \beta + 188370 \beta^{2} + 470400 \beta^{3} + 269325 \beta^{4} \\
\\
\theta_{6}^{7}(\beta) &= 1 + 114 \beta + 1452 \beta^{2} + 4400 \beta^{3} + 3708 \beta^{4} 
+ 720 \beta^{5} \\
\theta_{6}^{8}(\beta) &= 8 + 1031 \beta + 14976 \beta^{2} + 52756 \beta^{3} + 53296 \beta^{4} 
+ 13068 \beta^{5} \\
\theta_{6}^{9}(\beta) &= 36 + 5175 \beta + 84430 \beta^{2} + 339000 \beta^{3} + 399180 \beta^{4} 
+ 118124 \beta^{5} \\
\theta_{6}^{10}(\beta) &= 120 + 19035 \beta + 344620 \beta^{2} + 1553430 \beta^{3} + 2088840 \beta^{4} 
+ 723680 \beta^{5} \\
\\
\theta_{7}^{8}(\beta) &= 1 + 240 \beta + 5610 \beta^{2} + 32120 \beta^{3} + 58140 \beta^{3} + 33984 \beta^{5} + 
5040 \beta^{6} \\
\theta_{7}^{9}(\beta) &= 9 + 2406 \beta + 62910 \beta^{2} + 407880 \beta^{3} + 852180 \beta^{4} + 592056 \beta^{5} 
+ 109584 \beta^{6} \\
\theta_{7}^{10}(\beta) &= 45 + 13260 \beta + 383475 \beta^{2} + 2777040 \beta^{3} + 6577080 \beta^{4} 
+ 5292600 \beta^{5} + 1172700 \beta^{6} \\
\\
\theta_{8}^{9}(\beta) &= 1 + 494 \beta + 19950 \beta^{2} + 195800 \beta^{3} + 644020 \beta^{4} + 785304 \beta^{5}
+ 341136 \beta^{6} + 40320 \beta^{7} 
\end{align*}
\begin{align*}
\theta_{8}^{10}(\beta) &= 10 + 5441 \beta + 242364 \beta^{2} + 2645150 \beta^{3} + 9799480 \beta^{4} 
+ 13711620 \beta^{5} \\
& \hspace{5mm} + 7028784 \beta^{6} + 1026576 \beta^{7} \\
\\
\theta_{9}^{10}(\beta) &= 1 + 1004 \beta + 67260 \beta^{2} + 1062500 \beta^{3} + 5765500 \beta^{4} 
+ 12440064 \beta^{5} \\
& \hspace{5mm} + 11026296 \beta^{6} + 3733920 \beta^{7} + 362880 \beta^{8} \\
\end{align*}

\hspace{5mm} 
What is first noticed is that the $\beta^{0}$ coefficients of the $\theta_{k}^{r}(\beta)$
polynomials have the form $\binom{r}{k}$. These correspond to sequences A000217, A000292, A000332, A000389, A000579,
A000580, A000581, A000582, A001287, respectively for $k = 2, \cdots, 10$. By comparing the coefficients of $\theta_{2
}^{r}(\beta)$ to sequences A000292 and A000914 it is determined that $\theta_{2}^{r}(\beta) = \binom{r}{3} + s(r+1,r-1)
\, \beta = \binom{r}{3} \, \left( 1 + \frac{1}{4} \, (3r - 1) \beta \right)$, where $s(n,m)$ are Stirling numbers of
the first kind. Comparing the coefficients of $\theta_{3}^{r}(\beta)$ to A000332 and A001303 leads to
\begin{align}\label{e11} 
\theta_{3}^{r}(\beta) = \binom{r}{4} + 2 r \, \binom{r}{4} \, \beta + s(r+2, r-1) \, \beta^2 = 
\binom{r}{4} \, \left( 1 + 2 r \beta + \frac{r(r-1)}{2} \, \beta^2 \right).
\end{align}
Similarly comparing to A000389 and A000915 $\theta_{4}^{r}(\beta)$ has the form 
\begin{align}\label{e12} 
\theta_{4}^{r}(\beta) &= \binom{r}{5} + \alpha_{1}^{4} \, \beta + \alpha_{2}^{4} \beta^2 + s(n+3, r-1) \, 
\beta^3 \nonumber\\
&= \binom{r}{5} \, \left( 1 + \frac{1}{6} \, (25 r + 7) \beta + \frac{1}{6} \, (15 r^2 - 5 r - 2) \, \beta^2 
+ \frac{1}{48} \, ( 15 r^3 - 30 r^2 + 5 r + 2) \, \beta^3 \right). 
\end{align}
Comparing A000579 and A053567 then
\begin{align}\label{e13} 
\theta_{5}^{r}(\beta) &= \binom{r}{6} \left( 1 + (8 r + 4) \, \beta + \frac{1}{4} \, (35 r^2 + 9 r - 2) \, \beta^2
+ \frac{1}{2} \, r (5 r^2 - 5 r - 2) \, \beta^3 \right. \nonumber\\ 
& \hspace{15mm} \left. + \frac{1}{16} \, r(r-1)(3 r^2 - 7 r - 2) \, \beta^3 \right).
\end{align}
Further exact forms of the polynomials $\theta_{k}^{r}(\beta)$ can be obtained by expanding the list of polynomials.
This is currently beyond the scope of this work. 

\vspace{5mm}

\section{Developement of \texorpdfstring{S(r, $nx + r \beta$, $\beta$)}{Lg}}
Since $S(r, \alpha, \beta)$ has been desired it can be quickly determined that a similar form will be obtained
when $\alpha \to y + r \beta$, where $y = n x$. The listing of polynomials will be given in the following list. 
The polynomials are given by:
\begin{align}\label{e14} 
S(0, y, \beta) &= \frac{1}{y} \nonumber \\
S(1, y + \beta, \beta) &= p = \frac{1}{1-\beta} \nonumber\\
S(r, y + r \beta, \beta) &= p^{r} \left[ y^{r-1} + \sum_{k=1}^{r-1} \phi_{k}^{r}(\beta) \, y^{r-k-1}
\, \beta^{k} \, p^{k} \right], \hspace{5mm} r \geq 2.
\end{align}
The list of polynomials in $\beta$ is:

\begin{align*}
\begin{tabular}{c|ccc}
\multicolumn{1}{c}{$r$} & \multicolumn{1}{c}{$\phi_{1}^{r}(\beta)$} &
\multicolumn{1}{c}{$\phi_{2}^{r}(\beta)$} & \multicolumn{1}{c}{$\phi_{3}^{r}(\beta)$} \\
\hline
2 & 1 \, (2 - $\beta$) & 0 & 0 \\
3 & 3 \, (2 - $\beta$) & 9 - 8 $\beta$ + 2 $\beta^{2}$ & 0 \\
4 & 6 \, (2 - $\beta$) & 48 - 44 $\beta$ + 11 $\beta^{2}$ & 64 - 79 $\beta$ + 36 $\beta^{2}$ - 6 $\beta^{3}$ \\
5 & 10 \, (2 - $\beta$) & 150 - 140 $\beta$ + 35 $\beta^{2}$ & 500 - 645 $\beta$ + 300 $\beta^{2}$ - 50 $\beta^{3}$ \\
6 & 15 \, (2 - $\beta$) & 360 - 340 $\beta$ + 85 $\beta^{2}$ & 2160 - 2865 $\beta$ + 1350 $\beta^{2}$ - 225 $\beta^{3}$ \\
7 & 21 \, (2 - $\beta$) & 735 - 700 $\beta$ + 175 $\beta^{2}$ & 6860 - 9275 $\beta$ + 4410 $\beta^{2}$ - 735 $\beta^{3}$ \\
8 & 28 \, (2 - $\beta$) & 1344 - 1288 $\beta$ + 322 $\beta^{2}$ & 17920 - 24570 $\beta$ + 11760 $\beta^{2}$ - 1960 $\beta^{3}$ \\
9 & 36 \, (2 - $\beta$) & 2268 - 2184 $\beta$ + 546 $\beta^{2}$ & 40824 - 56574 $\beta$ + 27216 $\beta^{2}$ - 4536 $\beta^{3}$ \\
10 & 45 \, (2 - $\beta$) & 3600 - 3480 $\beta$ + 870 $\beta^{2}$ & 84000 - 117390 $\beta$ + 56112 $\beta^{2}$ - 13272 $\beta^{3}$ \\
\end{tabular}
\end{align*}

\hspace{5mm}
\begin{align*}
\phi_{4}^{5}(\beta) &= 625 - 974 \beta + 622 \beta^{2} - 192 \beta^{3} + 24 \beta^{4} \\
\phi_{4}^{6}(\beta) &= 6480 - 10614 \beta + 6997 \beta^{2} - 2192 \beta^{3} + 274 \beta^{4}  \\
\phi_{4}^{7}(\beta) &= 36015 - 60984 \beta + 41062 \beta^{2} - 12992 \beta^{3} + 1624 \beta^{4}  \\
\phi_{4}^{8}(\beta) &= 143360 - 248584 \beta + 169932 \beta^{2} - 54152 \beta^{3} + 6769 \beta^{4}  \\
\phi_{4}^{9}(\beta) &= 459270 - 810684 \beta + 560532 \beta^{2} - 179592 \beta^{3} + 22449 \beta^{4} \\
\phi_{4}^{10}(\beta) &= 1260000 - 2255148 \beta + 1537914 \beta^{2} - 700224 \beta^{3} + 292593 \beta^{4} \\
\\
\phi_{5}^{6}(\beta) &= 7776 - 14543 \beta + 11758 \beta^{2} - 5126 \beta^{3} + 1200 \beta^{4} 
- 120 \beta^{5} \\
\phi_{5}^{7}(\beta) &= 100842 - 198639 \beta + 166824 \beta^{2} - 74508 \beta^{3} + 17640 \beta^{3} 
- 1764 \beta^{5} \\
\phi_{5}^{8}(\beta) &= 688128 - 1405348 \beta + 1211728 \beta^{2} - 550326 \beta^{3} + 131320 \beta^{4} 
- 13132 \beta^{5} \\
\phi_{5}^{9}(\beta) &= 3306744 - 6935460 \beta + 6096720 \beta^{2} - 2803290 \beta^{3} + 672840 \beta^{4} 
- 67284 \beta^{5} \\
\phi_{5}^{10}(\beta) &= 12600000 - 26972190 \beta + 24064740 \beta^{2} - 11170530 \beta^{3} + 2693250 \beta^{4}
- 269325 \beta^{5} \\
\\
\phi_{6}^{7}(\beta) &= 117649 - 255828 \beta + 248250 \beta^{2} - 137512 \beta^{3} + 45756 \beta^{3} - 8640 \beta^{5} 
+ 720 \beta^{6} \\
\phi_{6}^{8}(\beta) &= 1835008 - 4203576 \beta + 4251207 \beta^{2} - 2426848 \beta^{3} + 823092 \beta^{4} 
- 156816 \beta^{5} + 13068 \beta^{6} \\
\phi_{6}^{9}(\beta) &= 14880348 - 35392104 \beta + 36865845 \beta^{2} - 21501260 \beta^{3} + 7392480 \beta^{4} 
- 1417488 \beta^{5} \nonumber\\
& \hspace{10mm} + 118124 \beta^{6} \\
\phi_{6}^{10}(\beta) &= 84000000 - 205539480 \beta + 218906235 \beta^{2} - 129746180 \beta^{3} + 45069630 \beta^{4}
\nonumber\\
& \hspace{10mm} - 8684160 \beta^{5} + 723680 \beta^{6} \\
\\
\phi_{7}^{8}(\beta) &= 2097152 - 5187775 \beta + 5846760 \beta^{2} - 3892430 \beta^{3} + 1651480 \beta^{4} 
- 445572 \beta^{5} + 70560 \beta^{6} \nonumber\\
& \hspace{10mm} - 5040 \beta^{7} \\
\theta_{7}^{9}(\beta) &= 38263752 - 99640521 \beta + 117218964 \beta^{2} - 80719830 \beta^{3} + 35097840 \beta^{4} 
- 9617772 \beta^{5} \\
& \hspace{5mm} + 1534176 \beta^{6} - 109584 \beta^{7} 
\end{align*}
\begin{align*}
\phi_{7}^{10}(\beta) &= 360000000 - 973850355 \beta + 1182292710 \beta^{2} - 834468975 \beta^{3} + 369371340 \beta^{4}
\nonumber\\
& \hspace{10mm} - 102373620 \beta^{5} + 16417800 \beta^{6} - 1172700 \beta^{7} \\
\\
\phi_{8}^{9}(\beta) &= 43046721 - 119214746 \beta + 152606870 \beta^{2} - 118016760 \beta^{3} + 60289700 \beta^{4} 
- 20808776 \beta^{5} \\
& \hspace{5mm} + 4728816 \beta^{6} - 645120 \beta^{7} + 40320 \beta^{8} \\
\phi_{8}^{10}(\beta) &= 900000000 - 2620523090 \beta + 3503248241 \beta^{2} - 2808432336 \beta^{3} 
+ 1475733950 \beta^{4} \nonumber\\
& \hspace{10mm} - 519839720 \beta^{5} + 119671020 \beta^{6} - 16425216 \beta^{7} 
+ 1026576 \beta^{8} \\
\\
\phi_{9}^{10}(\beta) &= 1000000000 - 3062575399 \beta + 3793195314 \beta^{2} - 4159084510 \beta^{3} 
+ 16269598360 \beta^{4} \nonumber\\
& \hspace{1mm} - 46433648700 \beta^{5} + 68381809464 \beta^{6} - 54285272064 \beta^{7} 
+ 22150947840 \beta^{8} - 3620510880 \beta^{9}
\end{align*}

\hspace{5mm} 
It is quickly determined that $\phi_{1}^{r}(\beta) = \binom{r}{2} \, (2 - \beta)$. The next few exact forms 
are given by:
\begin{align}
\phi_{2}^{r}(\beta) &= \frac{1}{4} \, \binom{r}{3} \, (12 r - 4 (3 r -1) \beta + (3 r -1) \beta^{2} ) 
\label{e15} \\ 
\phi_{3}^{r}(\beta) &= \frac{1}{2} \, \binom{r}{4} \, ( 8 r^{2} - 2(6 r^{2} - 4 r -1) \beta + 6 r (r-1) \beta^{2}
- r(r-1) \beta^{3} ) \label{e16} \\ 
\phi_{4}^{r}(\beta) &= \frac{1}{48} \, \binom{r}{5} \, ( 240 r^{3} - 48 (10 r^{3} - 10 r^{2} - 5 r -1) \beta
+ 8(45 r^{3} - 75 r^{2} - 5 r + 7) \beta^{2} \nonumber\\
& \hspace{5mm} - 8 (15 r^{3} - 30 r^{2} + 5 r + 2) \beta^{3} + (15 r^{3} 
- 30 r^{2} + 5 r + 2) \beta^{4} ) \label{e17} \\ 
\phi_{5}^{r}(\beta) &= \frac{1}{16} \, \binom{r}{6} \, ( 96 r^{4} - 16 (15 r^{4} - 20 r^{3} - 15 r^{2} - 6 r - 1)
\beta + 16 ( 15 r^{4} - 35 r^{3} - 5 r^{2} + 9 r + 4 ) \beta^{2} \nonumber\\
& \hspace{5mm} + 4 (30 r^{4} - 90 r^{3} + 25 r^{2} + 27 r 
+ 2) \beta^{3} + 10 r (r-1) (3 r^{2} - 7 r -2) \beta^{4} \nonumber\\
& \hspace{10mm} + r (r-1) (3 r^{2} - 7 r -2) \beta^{5} ). \label{e18} 
\end{align}
Equations (\ref{e14} - \ref{e18}) will be essential in providing the exact forms of the 
first few Jain moments.

\vspace{5mm}
\section{Jain Moments \texorpdfstring{$B_{n}^{\beta}(t^{m}, x)$}{Lg}}
In the introduction section the first couple of moments were shown how they connect to $S(r, nx + r \beta, 
\beta)$. The general form is given by, where $y = n x$,
\begin{align}\label{e19} 
B_{n}^{\beta}(t^{m}, x) = \frac{y}{n^{m}} \, \sum_{r=1}^{m} S(m,r) \, S(r, y + r \beta, \beta), 
\end{align}
where $S(m,r)$ are the Stirling numbers of the second kind. Up to $m=5$ the moments are
\begin{align*}
B_{n}^{\beta}(t^{0},x) &= 1 \\
B_{n}^{\beta}(t^{1},x) &= \frac{y p}{n} \\
B_{n}^{\beta}(t^{2},x) &= \frac{y p^{2}}{n^{2}} \, ( y + p ) \\
B_{n}^{\beta}(t^{3},x) &= \frac{y p^{3}}{n^{3}} \, ( y^{2} + 3 y p + (1 + 2 \beta) p^{2} ) 
\end{align*}
\begin{align*}
B_{n}^{\beta}(t^{4},x) &= \frac{y p^{4}}{n^{4}} \, (y^{3} + 6 y^{2} p + (7 + 8 \beta) \, y p^{2} + (1 + 8 \beta + 
6 \beta^{2}) p^{3} ) \\
B_{n}^{\beta}(t^{5},x) &= \frac{y p^{5}}{n^{5}} \, (y^{4} + 10 y^{3} p + 5 (5 + 4 \beta) \, y^{2} p^{2} + 15 (1 + 
4 \beta + 2 \beta^{2}) \, y p^{3} \nonumber\\
& \hspace{10mm} + (1 + 22 \beta + 58 \beta^{2} + 24 \beta^{3}) p^{4} ).  
\end{align*}
From this list it can be seen that the general form is
\begin{align}\label{e20} 
B_{n}^{\beta}(t^{m}, x) = \frac{y \, p^{m}}{n^{m}} \, \left[ y^{m-1} + \binom{m}{2} \, y^{m-2} \, p + 
\sum_{k=2}^{m-2} \sigma_{k}^{m}(\beta) \, y^{m-k-1} \, p^{k} + B_{m-1}(\beta) \, p^{m-1} \right],
\end{align}
where $B_{n}(x)$ are the Eulerian polynomials of the second kind. The coefficients of the Eulerian 
polynomials of the second kind are given by sequence A008517. The first few $\sigma_{k}^{m}(\beta)$ 
polynomials are given by
\begin{align}
\sigma_{2}^{m}(\beta) &= \frac{1}{4} \, \binom{m}{3} \, ( 3 m - 5 + 8 \beta) \label{e21} \\ 
\sigma_{3}^{m}(\beta) &= \frac{1}{2} \, \binom{m}{4} \, ( (m-2)(m-3) + 8(m-2)\beta + 12 \beta^{2}) \label{e22} \\
\sigma_{4}^{m}(\beta) &= \frac{1}{48} \, \binom{m}{5} \, ( (15 m^{3} - 150 m^{2} + 485 m - 502) + 16 (15 m^{2}
- 95 m + 116) \beta \nonumber\\ 
& \hspace{10mm} + 16 (65 m - 151) \beta^{2} + 1152 \beta^{3} ) \label{e23} \\ 
\sigma_{5}^{m}(\beta) &= \frac{1}{16} \, \binom{m}{6} \, ( (3 m^{4} - 50 m^{3} + 305 m^{2} - 802 m + 760) 
+ 16 (5 m^{3} - 55 m^{2} \nonumber\\ 
& \hspace{5mm} + 196 m - 224) \beta + 8 (85 m^{2} - 537 m + 818) \beta^{2} + 192 ( 11 m - 29) \beta^{3}
+ 120 \beta^{4} ). \label{e24} 
\end{align}
With the exact forms of (\ref{e20} - \ref{e24}) many of the polynomial coefficients in the 
moments can be quickly determined. The remaining terms not presented here, in exact form, required calculations
beyond the scope of this work. While most of the articles written involving the use of Jain moments use up to
$m = 5$ the above work yields all necessary terms. The listing of moments up to $m=10$ is:

\begin{align*}
B_{n}^{\beta}(t^{0},x) &= 1 \\
B_{n}^{\beta}(t^{1},x) &= \frac{y p}{n} \\
B_{n}^{\beta}(t^{2},x) &= \frac{y p^{2}}{n^{2}} \, ( y + p ) \\
B_{n}^{\beta}(t^{3},x) &= \frac{y p^{3}}{n^{3}} \, ( y^{2} + 3 y p + (1 + 2 \beta) p^{2} ) \\
B_{n}^{\beta}(t^{4},x) &= \frac{y p^{4}}{n^{4}} \, (y^{3} + 6 y^{2} p + (7 + 8 \beta) \, y p^{2} + (1 + 8 \beta + 
6 \beta^{2}) p^{3} ) \\
B_{n}^{\beta}(t^{5},x) &= \frac{y p^{5}}{n^{5}} \, (y^{4} + 10 y^{3} p + 5 (5 + 4 \beta) \, y^{2} p^{2} + 15 (1 + 
4 \beta + 2 \beta^{2}) \, y p^{3} \nonumber\\
& \hspace{10mm} + (1 + 22 \beta + 58 \beta^{2} + 24 \beta^{3}) p^{4} ) 
\end{align*}

\newpage
\begin{align*}
B_{n}^{\beta}(t^{6}, x) &= \frac{y p^{6}}{n^{6}} \, ( y^{5} + 15 y^{4} p + 5(13 + 8 \beta) y^{3} p^{2} + 
30(3 + 8 \beta + 3 \beta^{2}) y^{2} p^{3} + (31 + 292 \beta \nonumber\\
& \hspace{5mm} + 478 \beta^{2} + 144 \beta^{3}) y p^{4} + 
(1 + 52 \beta + 328 \beta^{2} + 444 \beta^{3} + 120 \beta^{4}) p^{5} ) \nonumber\\
B_{n}^{\beta}(t^{7}, x) &= \frac{y p^{7}}{n^{7}} \, ( y^{6} + 21 y^{5} p + 70(2 + \beta) y^{4} p^{2} +
70(5 + 10 \beta + 3 \beta^{2}) y^{3} p^{3} + 7(43 + 256 \beta \nonumber\\
& \hspace{5mm} + 304 \beta^{2} + 72 \beta^{3}) y^{2} p^{4} 
+ 21 ( 3 + 56 \beta + 204 \beta^{2} + 192 \beta^{3} + 40 \beta^{4}) y p^{5} + (1 \nonumber\\
& \hspace{10mm} + 114 \beta + 1452 \beta^{2} + 4400 \beta^{3} + 3708 \beta^{4} + 720 \beta^{5}) p^{6} ) \\
B_{n}^{\beta}(t^{8}, x) &= \frac{y p^{8}}{n^{8}} \, ( y^{7} + 28 y^{6} p + 14(19 + 8 \beta) y^{5} p^{2} 
+ 210(5 + 8 \beta + 2 \beta^{2}) y^{4} p^{3} + 21(81 + 352 \beta \nonumber\\
& \hspace{5mm} + 328 \beta^{2} + 64 \beta^{3}) y^{3} p^{4} 
+ 42(23 + 256 \beta + 654 \beta^{2} + 472 \beta^{3} + 80 \beta^{4}) y^{2} p^{5} +  (127 \nonumber\\
& \hspace{10mm} + 4272 \beta + 29172 \beta^{2} + 58832 \beta^{3} + 36972 \beta^{4} + 5760 \beta^{5}) y p^{6} 
+ (1 + 240 \beta \nonumber\\ 
& \hspace{10mm} + 5610 \beta^{2} + 32120 \beta^{3} + 58140 \beta^{4} + 33984 \beta^{5} + 5040 \beta^{6}) p^{7} ) \\
B_{n}^{\beta}(t^{9}, x) &= \frac{y p^{9}}{n^{9}} \, ( y^{8} + 36 y^{7} p + 42(11 + 4 \beta) y^{6} p^{2} 
+ 126(21+ 28 \beta + 6 \beta^{2}) y^{5} p^{3} + 21(331 \nonumber\\
& \hspace{5mm} + 1132 \beta + 868 \beta^{2} + 144 \beta^{3}) y^{4} p^{4} 
+ 210(37 + 292 \beta + 574 \beta^{2} + 336 \beta^{3} + 48 \beta^{4}) y^{3} p^{5} \nonumber\\
& \hspace{5mm} + 5(605 + 11388 \beta + 52656 \beta^{2} + 79504 \beta^{3} + 39852 \beta^{4} + 5184 \beta^{5})
y^{2} p^{6} + 15(17 \nonumber\\ 
& \hspace{5mm} + 972 \beta + 11286 \beta^{2} + 41200 \beta^{3} + 53964 \beta^{4} + 24672 \beta^{5} + 3024 \beta^{6}) 
y p^{7} + (1 + 494 \beta \nonumber\\
& \hspace{5mm} + 19950 \beta^{2} + 195800 \beta^{3} + 644020 \beta^{4} + 785304 \beta^{5} + 341136 
\beta^{6} + 40320 \beta^{7}) p^{8} ) \\
B_{n}^{\beta}(t^{10}, x) &= \frac{y p^{10}}{n^{10}} \, ( y^{9} + 45 y^{8} p + 30(25 + 8 \beta) y^{7} p^{2} 
+ 420(14 + 16 \beta + 3 \beta^{2}) y^{6} p^{3} + 21(1087 \nonumber\\
& \hspace{5mm} + 3064 \beta + 1996 \beta^{2} + 288 \beta^{3}) y^{5} p^{4}
+ 315(135 + 824 \beta + 1316 \beta^{2} + 648 \beta^{3} + 80 \beta^{4}) y^{4} p^{5} \nonumber\\
& \hspace{5mm} + 5(6821 + 88248 \beta + 307680 \beta^{2} + 371152 \beta^{3} + 154764 \beta^{4} + 17280 
\beta^{5}) y^{3} p^{6} \nonumber\\
& \hspace{5mm} + 30(311 + 9288 \beta + 70401 \beta^{2} + 187864 \beta^{3} + 192732 \beta^{4} + 72384 \beta^{5}
+ 7560 \beta^{6} )y^{2} p^{7} \nonumber\\
& \hspace{5mm} + (511 + 47804 \beta + 886650 \beta^{2} + 5317280 \beta^{3} + 12298780 \beta^{4} + 11470464 \beta^{5} 
\nonumber\\
& \hspace{5mm} + 4034736 \beta^{6} + 403200 \beta^{7} ) y p^{8} + (1 + 1004 \beta + 67260 \beta^{2} + 1062500 \beta^{3} 
\nonumber\\
& \hspace{5mm} + 5765500 \beta^{4} + 12440064 \beta^{5} + 11026296 \beta^{6} + 3733920 \beta^{7} 
+ 362880 \beta^{8}) p^{9} )
\end{align*}

\vspace{10mm}


\end{document}